\documentclass[a4paper,11pt]{amsart}

\usepackage{amssymb,enumerate,color}
\usepackage{color}
\usepackage{latexsym,amsbsy,bm}

\input epsf

\renewcommand{\O}{\mathcal{O}}

\newcommand {\SC} {{\mathbb C}}

\newcommand {\SN} {{\mathbb N}}
\newcommand {\SR} {{\mathbb R}}
\newcommand {\SRd} {{\SR^d}}

\newcommand {\SX} {{\mathbb X}}

\renewcommand {\phi} {{\varphi}}

\newcommand {\al} {{\alpha}}

\newcommand {\Dt} {{\Delta}}
\newcommand {\e} {{\varepsilon}}
\newcommand {\ga} {{\gamma}}
\newcommand {\Ga} {{\Gamma}}
\newcommand {\la} {{\lambda}}

\newcommand{\cD}{\mathcal D}

\newcommand {\tf} {{\tilde f}}

\newcommand {\re} {{\Re\!\mbox{\small\it e}\,}}

\def\span{\mathop{\rm span}}

\numberwithin{equation}{section}
\newtheorem{theorem}{Theorem}[section]
\newtheorem{lemma}[theorem]{Lemma}

\newtheorem{Remark}[theorem]{Remark}
\newtheorem{proposition}[theorem]{Proposition}

\newtheorem{example}[theorem]{Example}

\newtheorem{problem}{Problem}

\newcommand {\Proof} {\noindent{\bf P{\footnotesize\bf ROOF}: } \ }
\newcommand {\ProofEnd} {
             \begin{flushright} \vskip -0.2in $\Box$ \end{flushright}}

\newcommand{\bk}{{\bf k}}

\newcommand{\Ba}[1]{\begin{array}{#1}}
\newcommand{\Ea}{\end{array}}
\newcommand{\Be}{\begin{equation}}
\newcommand{\Ee}{\end{equation}}
\newcommand{\Bea}{\begin{eqnarray}}
\newcommand{\Eea}{\end{eqnarray}}
\newcommand{\Beas}{\begin{eqnarray*}}
\newcommand{\Eeas}{\end{eqnarray*}}
\newcommand{\Benu}{\begin{enumerate}}
\newcommand{\Eenu}{\end{enumerate}}
\newcommand{\Bi}{\begin{itemize}}
\newcommand{\Ei}{\end{itemize}}

\newcommand{\BR}{\begin{Remark} \em}
\newcommand{\ER}{\end{Remark}}
\newcommand{\BE}{\begin{example} \em}
\newcommand{\EE}{\end{example}}

\newcommand {\Ss} {\scriptstyle}

\newcommand {\Ds} {\displaystyle}

\renewcommand {\ae} {{\quad a.e.\,}}
\newcommand {\sae} {{\,a.e.\,}}
\newcommand {\mand} {{\quad\mbox{and}\quad}}
\renewcommand {\mid} {{\,\,\,\colon\,\,\,}}

\renewcommand{\th}{{\operatorname{th}\,}}

\newcommand{\sh}{{\operatorname{sh}\,}}

\renewcommand{\wp}{w^{-\frac1{p-1}}}

\newcommand{\Ml}{{\mathcal{M}^{\rm loc}}}
\newcommand{\bline}{{\bigskip

\noindent}}

\begin{document}

\title[2-weight problem for Poisson-Hermite]{A weak 2-weight problem for the Poisson-Hermite semigroup}
\author[Garrig\'os]{G. Garrig\'os}

\address{Gustavo Garrig\'os
\\
Departamento de Matem\'aticas
\\
Universidad de Murcia
\\
30100 Murcia, Spain. \emph{Email}: {\rm \texttt{gustavo.garrigos@um.es}}} 

\thanks{Author partially supported by grants MTM2010-16518, MTM2011-25377 and MTM2013-40945-P (Spain).}

\date{\today}


 \maketitle

\begin{abstract}
This survey is a slightly extended version of the lecture given by the author
at the \emph{VI International Course of Mathematical Analysis in Andaluc\'\i a} (CIDAMA), in 
September 2014. Most results form part of the paper \cite{GHSTV}, written jointly with S. Hartzstein, T. Signes, J.L. Torrea and B. Viviani.

\end{abstract}

\section{Introduction}  \setcounter{equation}{0}\setcounter{theorem}{0}
\setcounter{problem}{0}

Consider the following integral identity
\Be
e^{-t\sqrt L}=\tfrac{t}{\sqrt{4\pi}}\,\int_0^\infty
 e^{-\frac{t^2}{4s}}\,e^{-sL}\,\frac{ds}{s^{3/2}},\quad t>0\label{subL}\Ee
valid for all real numbers $L>0$.
If we allow $L$ be the infinitesimal generator of a ``heat'' semigroup $\{e^{-sL}\}_{s>0}$ in $L^2(\SR^d)$,
then \eqref{subL} defines, using the terminology 
in Stein's book \cite[Chapter II.2]{steinLP}, a \emph{subordinated Poisson semigroup}. Moreover, for suitably ``good'' functions $f:\SR^d\to\SC$, 
the formal \emph{Poisson integral} $u(t,\cdot)=e^{-t\sqrt L}f$ 
solves the partial differential equation
\[u_{tt} = Lu,\quad (t,x)\in(0,\infty)\times\SR^d, \quad \mbox{ with } u(0)=f.\]
A relevant question is then to find, for each operator $L$, the most general class of functions $f$ for which 
the Poisson integrals $u(t,x)=e^{-t\sqrt L}f(x)$ satisfy
\Bi\item[(i)] $u(t,x)$ is well-defined and belongs to $C^\infty((0,\infty)\times\SR^d)$
\item[(ii)] $u(t,x)$ satisfies the pde $u_{tt}=Lu$ in $(0,\infty)\times\SR^d$
\item[(iii)]  There exists $\lim_{t\to0^+}u(t,x)=f(x)$, for a.e. $x\in\SR^d$.
\Ei 
In the classical setting, corresponding to the Laplace operator $L=-\Dt$ in $\SR^d$,
the largest class of admissible initial data $f$ is the weighted space 
\Be L^1(\phi)=\Big\{f\mid \int_{\SR^d}{|f(x)|}\phi(x)\,dx<\infty\Big\}, 
\label{pd1}\Ee
 with $\phi(x)={(1+|x|)^{-(d+1)}}$, and the assertions (i)-(iii) can easily be proved from the explicit form of the Poisson kernel.

For general operators $L$, however, the kernel will not be so explicit, 
and investigating such results requires very precise estimates of the subordinated integrals
 in \eqref{subL}, as well as of the associated maximal operators.

In this work we take up this question for a collection of Hermite operators in $\SR^d$
\Be
L=-\Dt +|x|^2 +m,  \quad \mbox{ with } m\geq -d.
\label{L}\Ee
We shall also consider a slightly more general family of partial differential equations:
\Be
u_{tt} \,+\,\tfrac{1-2\nu}t\,u_t\,=\,L u,\quad (t,x)\in(0,\infty)\times\SR^{d}, \quad \mbox{ with }\nu>0.
\label{pde}\Ee
The parameters $m$ and $\nu$ allow us to include various interesting cases, 
which can all be covered with essentially the same estimates. 
In particular, $m=0$ corresponds to the usual Hermite operator, while $m=-d$ leads to an $L$ which can be transformed\footnote{Note that $e^{{|x|^2}/2}L[e^{-{|x|^2}/2}u]=-\Dt u+2x\cdot\nabla u +(m+d)u$.} into the Ornstein-Uhlenbeck operator $-\Dt+2x\cdot\nabla$. Likewise, the parameter $\nu=1/2$ in \eqref{pde}  gives  the  usual Poisson equation, while for general $\nu$ it leads to a pde appearing in the theory of \emph{fractional laplacians}\footnote{The fractional operator $L^\nu$ can be recovered from \eqref{pde} and $u(0,x)=f(x)$  by the formula
$L^\nu f(x)=c_\nu\lim_{t\to0}t^{1-2\nu}u_t(t,x)$, at least for suitably good $f$;
see \cite[Thm 1.1]{StTo}.}.
  
\medskip

Our goal in this work is to give the most general conditions on a function $f:\SR^d\to\SC$ so that a meaningful solution to \eqref{pde} is given by the \emph{Poisson-like integral}
\Be
P_tf(x):=\tfrac{t^{2\nu}}{4^\nu\Gamma(\nu)}\,\int_0^\infty
 e^{-\frac{t^2}{4u}}\,\big[e^{-uL}f\big](x)\,\frac{du}{u^{1+\nu}},\quad t>0.
\label{poisf}\Ee
This subordinated integral is slightly more general than \eqref{subL}, and we justify its expression in $\S2$ below.
In our results, which we are about to state, the following function will play a crucial role
\Be
\phi(y)=\left\{\Ba{ll}
\Ds\frac {e^{-|y|^2/2}}{(1+|y|)^{\frac {d+m}2}[\ln(e+|y|)]^{1+\nu}} & \mbox{ if } m>-d\\
&\\
\Ds\frac {e^{-|y|^2/2}}{[\ln(e+|y|)]^{\nu}} & \mbox{ if } m=-d\\
\Ea\right.
\label{phi}\Ee

\

\begin{theorem}\label{th1}
For every $f\in L^1(\phi)$ the function $u(t,x)=P_tf(x)$ in \eqref{poisf} is defined by an absolutely convergent
integral such that
\Bi\item[(i)] $u(t,x)\in C^\infty((0,\infty)\times\SR^d)$ 
\item[(ii)] $u(t,x)$ satisfies the pde \eqref{pde}
\item[(iii)]  For a.e. $x\in\SR^d$, it holds $\lim_{t\to0^+}u(t,x)=f(x)$.
\Ei
Conversely, if a function $f\geq0 $  is such that the integral in \eqref{poisf} is finite for some $(t,x)\in(0,\infty)\times\SR^d$, then $f$ must necessarily belong to $L^1(\phi)$.
\end{theorem}

In particular, a function with growth as large as $f(y)=e^{\frac{|y|^2}{2}}/[(1+|y|)^d\ln(e+|y|)]$  has nicely convergent Poisson integrals, for all $m\geq -d$ and $\nu>0$.  This is in contrast with the 
classical case $L=-\Dt$ for which only a mild sublinear growth is allowed; see \eqref{pd1}.
It also illustrates that $L^1(\phi)$ is strictly larger than the ``gaussian'' space $L^1(\SR^d,e^{-\frac{|y|^2}2}dy)$, which was the natural domain for Poisson integrals considered by Muckenhoupt in \cite{Muck1}
(in the special case $\nu=1/2$, $m=-d$ and $d=1$).

\medskip

Our second goal is to investigate the following local maximal operators
\Be
\label{max} 
P_a^*f(x):=\sup_{0<t<a}\big|P_tf(x)\big|\;, \quad\mbox{with $a>0$ fixed.}\Ee 
These operators arise naturally in the $\sae$-pointwise convergence of $P_tf(x)\to f(x)$
as $t\to0$. In fact, the natural strategy to prove such convergence for all $f$ in a Banach space $\SX$, is to establish first the result in a dense class,
and next prove the boundedness of $P_a^*$ from $\SX$ into $L^{p, \infty}(v)$ (or even better into $L^p(v)$) for some weight $v>0$. 
It turns out that we can prove Theorem \ref{th1} without appeal to such maximal operators,
but it still makes sense to consider the following 

\begin{problem} \label{Prob}{\bf A weak 2-weight problem for the operator $P^*_a$.}
Given $a>0$ and  $1<p<\infty$, characterize the weights $w(x)>0$ for which
there exists some other weight $v(x)>0$ such that\Be
P^*_a:L^p(w)\to L^p(v) \quad \mbox{boundedly}.\label{bded}\Ee
\end{problem}

We named the problem ``weak'' in contrast with the ``strong'' (and more difficult) question of characterizing all pairs of weights $(w,v)$ for which \eqref{bded} holds. Such weak 2-weight problems,
for various classical operators, were considered in the early 80s by Rubio de Francia \cite{RdF1}
and Carleson and Jones \cite{CJ}, who found explicit answers for the Hardy-Littlewood maximal operator and the Hilbert transform.

\

Our second main result in \cite{GHSTV} gives an answer to Problem \ref{Prob}.

\begin{theorem}\label{th2}
Let $1<p<\infty$ and $a>0$ be fixed. Then, for a weight $w(x)>0$ the condition
\Be
\big\|w^{-\frac1p}\,\phi\,\big\|_{L^{p'}(\SR^d)}<\infty \label{Dp}
\Ee
is equivalent to the existence of some other weight $v(x)>0$ such that \eqref{bded} holds.
\end{theorem}

Condition \eqref{Dp} is easily seen to be equivalent to $L^p(w)\subset L^1(\phi)$.
So, the necessity of \eqref{Dp} in Theorem \ref{th2} is a consequence of the last sentence in Theorem \ref{th1}.
Concerning the sufficiency of \eqref{Dp}, we first point out that, from Theorem \ref{th1} (iii) and abstract results
due to Nikishin, there always exists a weight $u(x)>0$ such that 
\Be
P^*_a:L^p(w)\to L^{p,\infty}(u) \quad \mbox{boundedly}.\label{wbded}\Ee
The main contribution of Theorem \ref{th2} is to show that the weak-space $L^{p,\infty}(u)$ in \eqref{wbded} can be replaced by the strong space $L^p(v)$ (with perhaps another weight $v$). This is the main difficulty in the 2-weight  Problem~\ref{Prob} described above, and requires additional estimates to those needed in Theorem \ref{th1}.

\medskip

A last question regards the explicit form of the weight $v(x)$, whose existence, under the condition \eqref{Dp}, is asserted in Theorem \ref{th2}.
In \cite{GHSTV} we used a non-constructive procedure which nevertheless provided a size estimate.
Namely, for every $\sigma<1$ a weight $v=v_{\sigma}$ can be chosen such that
\Be
\big\|v^{-\frac\sigma p}\,\phi\,\big\|_{L^{p'}(\SR^d)}<\infty \label{Dpe}.
\Ee
Notice that this is ``almost'' the same integrability condition that $w(x)$ satisfies.
Here we state a new result, which provides an explicit expression for $v(x)$, and recovers in particular
the property \eqref{Dpe}. We shall  use the following \emph{local Hardy-Littlewood maximal operator}
\Be
\Ml f(x)=\sup_{r>0}\frac1{|B_r|}\int_{B_r(x)}|f(y)|\,\chi_{\{|y|\leq 3\max(|x|,1)\}}\,dy.
\label{Ml}\Ee

\

\begin{theorem}\label{th3}
Let $1<p<\infty$ be fixed, and let $w(x)>0$ be a weight satisfying \eqref{Dp}.
Then a family of weights $v(x)$ such that \eqref{bded} holds for all $a>0$ is given explicitly by

\Be
v(x)=\Big[\mathcal{M}^{\rm loc}\big(w^{-\frac{p'}p}\,e^{-\frac{p'|y|^2}2}\big)(x)\Big]^{-\frac\al{p'/p}}\,e^{-\frac{p|x|^2}2}\, (1+|x|)^{-N}, 
\label{vx}\Ee
\

\noindent provided $\al>1$ and $N> N_0$, for some $N_0=N_0(\al,p,d,m,\nu)$. 
\end{theorem}

We finally remark that, via the elementary identity 
\[e^{{|x|^2}/2}L[e^{-{|\,\cdot\,|^2}/2}u]\,=\,-\Dt u+2x\cdot\nabla u +(m+d)u\,=:\O,\] 
all the results in this paper admit corresponding statements with $L$ replaced by the
Ornstein-Uhlenbek type operator $\O$. These essentially amount to replace the exponentials
$e^{-|y|^2/2}$ (as in \eqref{phi} or \eqref{vx}) by the gaussians $e^{-|y|^2}$. 
We leave the simple verification to the interested reader.

\

The proof of Theorems \ref{th1} and \ref{th2} was given in \cite{GHSTV}, but we outline
the main steps below. Namely, in $\S2$ we justify why the integral formula in \eqref{poisf}
gives a solution to the pde \eqref{pde}. In $\S3$ we state the optimal kernel estimates 
which are behind these theorems, and outline the proof of Theorem~\ref{th1}, slightly modified with respect to \cite{GHSTV}.
The new results appear in $\S4$, where we solve a weak 2-weight problem for $\Ml$, and
present the proof of Theorem~\ref{th3} (which in turn implies Theorem~\ref{th2}).

\section{The subordinated integral}  \setcounter{equation}{0}\setcounter{theorem}{0}
\setcounter{problem}{0}

For $\nu\in\SR$, consider the following real-valued function
\Be
F_\nu(z):=\int_0^\infty e^{-u-\frac{z^2}{4u}}\,u^{\nu-1}\,du,\quad z>0,
\label{Fnu}\Ee
where the integral is absolutely convergent (actually for all $\nu\in\SC$ and $\re(z^2)>0$).
This integral is well-known in the theory of special functions, as it appears in the definition of the so-called modified Bessel function of the third kind $K_\nu(z)$. Namely, they are related by
\Be
F_\nu(z)=2\,(z/2)^\nu\, K_\nu(z);\label{Knu}
\Ee
see e.g. \cite[p. 183]{wat} or \cite[p. 119]{leb}.
In particular, $F_\nu$ satisfies the ordinary differential equation
\Be\label{ode}
F_\nu''(z)+\frac{1-2\nu}z F_\nu'(z) \,=\, F_\nu(z).\Ee
We give next the elementary proof of \eqref{ode}, which does not depend on
the properties of $K_\nu$. Integrating by parts in \eqref{Fnu} we can write
\[
F_\nu(z) =  \int_0^\infty e^{-u} \Big(e^{-\frac{z^2}{4u}}\,u^{\nu-1}\Big)'\,du\,=\,
\int_0^\infty e^{-u} \big(\tfrac{z^2}{4u^2} +\tfrac{\nu-1}u\big)\,
e^{-\frac{z^2}{4u}}\,u^{\nu-1}\,du.\]
Taking derivatives inside the integral in \eqref{Fnu} we also have
\[
F'_\nu(z)= \int_0^\infty e^{-u-\frac{z^2}{4u}}\,\big(-\tfrac{z}{2u}\big)\,u^{\nu-1}\,du
\mand 
F''_\nu(z)= \int_0^\infty e^{-u-\frac{z^2}{4u}}\,\big(\tfrac{z^2}{4u^2}-\tfrac{1}{2u}\big)\,u^{\nu-1}\,du.
\]
From these identities \eqref{ode} follows easily. 
Moreover we have the following

\begin{lemma}
Let $\nu$ and $L$ be positive real numbers. Then, the function 
$u(t)=\tfrac1{\Ga(\nu)}F_\nu(t\sqrt L)$, $t>0$, satisfies the differential equation
\Be
u''(t) + \frac{1-2\nu}t\,u'(t)\,=\,Lu(t),\quad \mbox{with }\lim_{t\to0^+}u(t)=1. \quad 
\label{ode2}\Ee Moreover, the function $u(t)$ can also be written as 
\Be
u(t)= \frac{(t/2)^\nu}{\Ga(\nu)}\,\int_0^\infty 
e^{-\frac{t^2}{4v}-Lv}\,\frac{dv}{v^{1+\nu}}, \quad t>0.
\label{uFt}\Ee\end{lemma}
\Proof
If $\nu>0$, from \eqref{Fnu} and dominated convergence it follows that 
\[
\lim_{z\to0^+}F_\nu(z) =\int_0^\infty e^{-u}u^{\nu-1}\,du\,=\,\Ga(\nu).
\] 
It is then straightforward to derive \eqref{ode2} from this observation and \eqref{ode}.
To obtain the integral expression in \eqref{uFt}, first set $z^2=t^2L$ in \eqref{Fnu},
and then change variables $v=\frac{t^2}{4u}$.\ProofEnd

When $L$ is a positive self-adjoint differential operator which generates a semigroup 
$\{e^{-sL}\}_{s>0}$ in $L^2(\SRd)$, we may then consider the function\[
u(t,x)=\tfrac{(t/2)^{2\nu}}{\Gamma(\nu)}\,\int_0^\infty
 e^{-\frac{t^2}{4v}}\,\big[e^{-vL}f\big](x)\,\frac{dv}{v^{1+\nu}},\quad t>0.
\]
In view of \eqref{ode2}, this is a natural candidate to solve the pde \[
u_{tt} \,+\,\tfrac{1-2\nu}t\,u_t\,=\,L u,\quad (t,x)\in(0,\infty)\times\SR^{d}, 
\quad \mbox{ with }u(0,\cdot)=f\]
(and coincides with the definition we used in \eqref{poisf} for the Poisson integral $P_tf(x)$ associated with $L$). Theorem \ref{th1} will give a rigorous proof of this formal statement, at least for the Hermite operators in \eqref{L}. We refer to \cite{StTo} for more on this kind of arguments for general operators $L$.

\section{Estimates on the Poisson kernels}

Suppose $L$ is the infinitesimal generator of a semigroup of operators
$\{h_t=e^{-tL}\}_{t>0}$, in $L^2(\SRd)$, and that these are described by the integrals
\Be
h_tf(x)=\int_{\SR^d}h_t(x,y)f(y)\,dy,
\label{htf}\Ee
for suitable positive kernels $h_t(x,y)$. 
Then, for the family of subordinated operators $\{P_t\}_{t>0}$ defined in \eqref{poisf}, a formal computation leads to
\[
P_tf(x)=\int_{\SR^d}p_t(x,y)f(y)\,dy,
\]
with the corresponding  kernels given by the integrals
\Be
p_t(x,y)=\tfrac{(t/2)^{2\nu}}{\Gamma(\nu)}\,\int_0^\infty
 e^{-\frac{t^2}{4v}}\,h_v(x,y)\,\frac{dv}{v^{1+\nu}}.
\label{pt0}\Ee
If one is interested in \emph{optimal} estimates for such kernels $p_t(x,y)$,
two things become necessary: first, a precise \emph{a priori} knowledge of $h_v(x,y)$, and
next a careful analysis of the integrals \eqref{pt0}.

Such tasks are difficult to carry in full generality, so in this work we have considered the special case of the Hermite operators $L=-\Dt+|x|^2+m$, 
for which we can start with an \emph{explicit} expression of the associated heat kernels
\Be
h_v(x,y)=\,e^{-vm}\,[2\pi\sh{2v}]^{-\frac d2}\,{e^{-\frac{|x-y|^2}{2\th{2v}}-\th{v} \,x\cdot y}}\,,\quad v>0;
\label{mehler}\Ee
see e.g. \cite[(4.3.14)]{Than}\footnote{Note that $e^{-vL}=e^{-v(-\Dt+|x|^2)} e^{-vm}$, 
and the case $m=0$ corresponds to the usual Mehler kernel.}.  
To make this expression more manageable it is common to use the new variable $s=\th(v)$ (or equivalently $v=\frac12\log(\frac{1+s}{1-s})$), which after elementary computations allows us to write
\[
h_v(x,y)=\,\frac{(1-s)^{\frac{m+d}2}}{(1+s)^{\frac{m-d}2}}\, \frac{e^{-\frac14(\frac{|x-y|^2}s + s|x+y|^2)}}{(4\pi s)^{\frac d2}}\,.
\]
Inserting this into the integral \eqref{pt0} (with $dv=\frac{ds}
{1-s^2}$) one obtains
the expression
\Be
p_t(x,y)=\tfrac{(t/2)^{2\nu}}{(4\pi)^\frac{d}2\Ga(\nu)}\,\int_0^1
\frac{e^{-\frac{t^2}{2\ln\frac{1+s}{1-s}}}\,(1-s)^{\frac{m+d}2-1}\,e^{-\frac14(\frac{|x-y|^2}s+s|x+y|^2)}} {s^{\frac
d2}\,(1+s)^{\frac{m-d}2+1}\,\big(\frac12\ln\frac{1+s}{1-s}\big)^{1+\nu}}\,ds,\label{poi0}\Ee
This will be our starting formula for $p_t(x,y)$, from which we shall derive the necessary estimates
needed for Theorems \ref{th1}, \ref{th2} and \ref{th3}.  
These are summarized in the next two lemmas.

The first one gives, for fixed $t$ and $x$, the \emph{optimal} decay of $y\mapsto p_t(x,y)$ in terms of the function $\phi(y)$ defined in \eqref{phi}. We shall sketch its proof in $\S$\ref{Lem3.1} below.

\begin{lemma}\label{L3.1}
Given $t>0$ and $x\in\SR^d$, there exist $c_1(t,x)>0$ and $c_2(t,x)>0$ such
that \Be c_1(t,x)\, \phi(y) \,\leq\,p_t(x,y)\,\leq\,
c_2(t,x)\, \phi(y) \;,\quad \forall\;y\in\SR^d. \label{pk1}\Ee
\end{lemma}

The second lemma is a refinement of the upper bound in \eqref{pk1} with two main advantages: it is uniform in the variable $t$,  and it restricts to the ``local part'' the singularities of the kernel $p_t(x,y)$.
The proof of this more precise lemma is sketched in $\S$\ref{Lem3.2}.

\begin{lemma}
\label{L3.2} Given $x\in\SRd$,  the following estimate holds for all $t>0$ and $y\in\SR^d$
\Be\label{pk2} p_t(x,y)\leq\,
\frac{C_1(x)\,t^{2\nu}\,e^{-\frac{|y|^2}2}}{(t+|x-y|)^{d+2\nu}}\,\chi_{\big\{|y|\leq
3\max\{|x|,1\}\big\}}
\,+\,C_2(x)\,t^{2\nu}\,\phi(y),
\Ee
for some positive functions $C_1(x) \,\lesssim\,(1+|x|)^{2\nu+d-1}e^{\frac{|x|^2}2}$ and $C_2(x)\lesssim 1/\phi(x)$.
\end{lemma}

 Observe that, as a consequence of \eqref{pk2},
 we obtain the following bound for the maximal operators $P^*_a$
in \eqref{max}
\Be
 P^*_af(x) \lesssim  C_1(x)\,\Ml \big(fe^{-\frac{|y|^2}2}\big)(x) \,+\,
C_2(x)\,a^{2\nu}\,\int_{\SRd}|f(y)|\,\phi(y)\,dy,\label{P**}
\Ee
where $\Ml$ denotes the local Hardy-Littlewood maximal operator defined in \eqref{Ml}.

\subsection{Proof of Theorem \ref{th1}}\label{Sth1}
 Assuming the lemmas, we can sketch the proof of Theorem \ref{th1}.
First of all, it is a direct consequence of Lemma \ref{L3.1} that $P_t|f|(x)<\infty$  for some (or all) 
$t>0$ and $x\in\SRd$ if and only if $f\in L^1(\phi)$. This justifies that $f\in L^1(\phi)$
is the right setting for this problem. Observe also that taking derivatives with respect to $t$ in $p_t(x,y)$ 
slightly improves the decay of the kernel, and from here it is not difficult to deduce that (i) and (ii) must hold;
see the details in \cite[Proposition 4.4]{GHSTV}.

We shall be a bit more precise about the proof of (iii), that is the pointwise convergence 
\Be
\lim_{t\to0^+}P_tf(x)=f(x),\ae x\in\SRd
\label{ptwise}\Ee
for all $f\in L^1(\phi)$. We first claim that such convergence holds in the dense set 
$\cD=\span\{h_{\bk}\}_{\bk\in\SN^d}$, where $h_{\bk}(x)$ denote the $d$-dimensional Hermite functions 
(as in \cite[p.5]{Than}). These are eigenfunctions of $L=-\Dt+|x|^2+m$ with
\[
Lh_{\bk} \,=\,(2|\bk|+d+m)\,h_{\bk},\quad {\rm if}\;|\bk|=k_1+\ldots+k_d\geq0\,;
\]
see \cite[(1.1.28)]{Than}. Recall that the operators $h_t=e^{-tL}$ from the Hermite semigroup can be 
represented in two ways: as in \eqref{htf} with the Mehler kernel 
\eqref{mehler}, or equivalently  as 
\[
h_tf = \sum_{\bk\in\SN^d}e^{-(2|\bk|+d+m)t}\,\langle f,h_{\bk}\rangle \,h_{\bk},
\]
at least for $f\in\cD$; see \cite[(4.1.1)]{Than}. 
From this last formula and the results in $\S2$ one also deduces that
\Be
P_tf = \,\tfrac1{\Ga(\nu)}\,\sum_{\bk\in\SN^d}\,F_\nu\big(t\sqrt{2|\bk|+d+m}\big)\,\langle f,h_{\bk}\rangle \,h_{\bk},\quad f\in\cD.
\label{Ptsum}\Ee
This clearly implies \eqref{ptwise} when $f\in\cD$. 

To extend this convergence to all $f\in L^1(\phi)$ 
we shall argue as in \cite[Proposition 4.5]{GHSTV}. 
Namely, it suffices to show \eqref{ptwise} for $\sae |x|\leq R$, for every fixed $R\geq1$.
We split  $f\in L^1(\phi)$ by
\[
f=\,{f\chi_{\{|y|\leq 3R\}}\,+\,f\chi_{\{|y|>3R\}}}\,=\,f_0+f_1.\]
Using Lemma \ref{L3.2} we see that, for every $|x|\leq R$ \Bea \big| P_tf_1(x)\big| & \leq & 
\int_{|y|>3R}p_t(x,y)\,|f(y)|\,dy\nonumber\\
&\leq & C_{R}\,t^{2\nu}\,\int_{\SRd}|f(y)|\phi(y)\,dy\to0,\quad{\rm as}\;t\to0^+. \label{3lf2}\Eea 
On the other hand, Lemma \ref{L3.2} (or rather its consequence in \eqref{P**}) also implies that
\[
\sup_{0<t\leq1}|P_tf_0(x)|\,\leq\, C_R\,\Big(Mf_0(x) +\int |f_0|\phi\Big),\quad |x|\leq R,
\]
where $M$ denotes the usual Hardy-Littlewood maximal operator. Since
the right-hand side is a bounded operator from $L^1(B_{3R}(0))\to L^{1,\infty}(B_R(0))$,
a classical procedure\footnote{See e.g. \cite[Theorem 2.2]{duo}.} then gives, from the validity of \eqref{ptwise} in the dense class $\cD$, the existence of
$\lim_{t\to0^+}P_tf_0(x)=f(x)$ for $\sae|x|\leq R$.
This completes the proof of Theorem \ref{th1}.
\ProofEnd
\BR
Notice that the series representation of $P_tf$ in \eqref{Ptsum} allows us to reformulate \eqref{ptwise}
as a result on pointwise convergence of Hermite expansions
by a ``summability method'' (based on the function $F_\nu$ and the parameter $m$).
This is in the same spirit as the Poisson summability for Hermite expansions considered by Muckenhoupt 
in \cite{Muck1} (in the special case $\nu=1/2$, $m=-d$ and $d=1$).
Notice, however, that the integral representation of $P_tf(x)$ in \eqref{poisf} is much more versatile,
as it exists for functions in $f\in L^1(\phi)$ which may have $\langle f,h_{\bk}\rangle=\infty$ for all $\bk$.
\ER

\subsection{Proof of Lemma \ref{L3.1}}\label{Lem3.1}
For the sake of originality, we shall use a slightly different approach than the one given in \cite[Lemma 4.1]{GHSTV}.
In the Mehler kernel $h_v(x,y)$ we shall consider the ``more natural'' variable
$r=e^{-2v}$ (or equivalently $v=\frac12\ln\frac1r$), which leads to the formula\footnote{This change of variables is common in the Ornstein-Uhlenbeck setting; see e.g. \cite[(3.3)]{Muck1} or \cite{LS}.}
\[
h_v(x,y)=
\frac{r^{\frac{m+d}2}\,e^{-\frac{|x-ry|^2}{1-r^2}}}{\Ss[\pi (1-r^2)]^{d/2}}\;e^{\frac{|x|^2-|y|^2}2}\,.
\]
Inserting this into the integral defining $p_t(x,y)$ (with $dv=dr/(2r)$) one obtains
the expression
\Be
p_t(x,y)=\tfrac{t^{2\nu}}{2^\nu\pi^\frac{d}2\Ga(\nu)}\,\int_0^1
\frac{e^{-\frac{t^2}{2\ln\frac1r}}\,r^{\frac{m+d}2}\,
e^{-\frac{|x-ry|^2}{1-r^2}}} {\Ss(1-r^2)^{\frac
d2}\,\big(\ln\frac1r\big)^{1+\nu}}\,\tfrac{dr}r\;e^{\frac{|x|^2-|y|^2}2}.\label{pois}\Ee
Starting from this formula, we now argue as in \cite[Lemma 4.1]{GHSTV}.
We may assume that $|y|\geq 3\max\{1,|x|\}$ (since in the region $|y|\leq 3\max\{1,|x|\}$
one can bound $p_t(x,y)/\phi(y)$ above and below 
by positive functions of $t,x$). 

The main difficulty is to determine the values
of $r$ which carry the main contribution of the integral in \eqref{pois}. 
The leading term will be the exponential in the numerator, and as we shall see, 
it becomes largest when $r\approx |x|/|y|$. 

Consider first the region $1/2<r<1$. Since $|y|\geq 3|x|$, we have
\[
|ry-x|\geq |y|/2-|x|\geq |y|/6\quad \Longrightarrow\quad  e^{-\frac{|x-ry|^2}{1-r^2}}\,\leq\, e^{-\frac{|y|^2}{36}},\]
so the leading exponential becomes quite small in this part. Since  we also have\Be
\ln\frac1r\approx
\, 1-r,\quad r\in[1/2,1],
\label{expts}\Ee we can estimate the corresponding integral in \eqref{pois} by 
\Bea \int_{1/2}^1\ldots& \lesssim & 
\,e^{\frac{|x|^2-|y|^2}2}\,e^{-\frac{|y|^2}{36}}\, t^{2\nu}\,\int_{1/2}^1
\frac{e^{-\frac{ct^2}{1-r}}} {(1-r)^{\frac
{d}2+\nu+1}}\,dr\nonumber\\ \mbox{\tiny [$u=\frac{t^2}{1-r}$]} 
& \leq & \,e^{\frac{|x|^2-|y|^2}2}\,e^{-\frac{|y|^2}{36}}\,
t^{-d}\int_0^\infty e^{-cu}\,u^{\frac d2+\nu-1}\,du\nonumber\\
&\lesssim & e^{\frac{|x|^2}2}\,t^{-d}\,\phi(y)\,.\label{I0}\Eea
Suppose now that $0<r<\frac12$. As we shall see, the main contribution occurs here, precisely when $r\approx \frac{|x|}{|y|}$. We first consider the range $2\frac{|x|}{|y|}<r<1/2$, where we can estimate
\[
|x-ry|\geq r|y|-|x|\geq r|y|/2\quad \Longrightarrow \quad 
e^{-\frac{|x-ry|^2}{1-r^2}}\,\leq\, e^{-\frac{r^2|y|^2}{4}}.\]
Thus 
\Beas 
\int_{2\frac{|x|}{|y|}}^{1/2}\ldots\quad & \lesssim & \,e^{\frac{|x|^2-|y|^2}2}\, t^{2\nu}\,
\int_{2\frac{|x|}{|y|}}^{1/2}
\frac{r^{\frac{m+d}2}\,e^{-\frac{r^2|y|^2}4}}{\big(\ln\frac1r\big)^{1+\nu}}\,\frac{dr}r\\
\mbox{{\tiny $\big[u=r|y|\big]$}} & \leq & t^{2\nu}\,e^{\frac{|x|^2-|y|^2}2}\, |y|^{-\frac{m+d}2}\,
\int_{0}^{\frac{|y|}2}
\frac{u^{\frac{m+d}2}\,e^{-u^2/4}}{\big(\ln\frac{|y|}{u}\big)^{1+\nu}}\,\frac{du}u.
\Eeas
Note that when $m+d>0$ the last integral\footnote{When 
$m=-d$, the integral still converges and is controlled by $1/(\ln|y|)^{\nu}$.} has its major contribution at $u\approx1$, so
it can be estimated by $1/(\ln|y|)^{1+\nu}$.  
Thus, overall one obtains
\Be \int_{2\frac{|x|}{|y|}}^{1/2}\ldots\,\lesssim\,
\,t^{2\nu}\,e^{\frac{|x|^2}2}\,\phi(y)\,.\label{I1a}\Ee
Finally, in the range $0<r<2\frac{|x|}{|y|}$ we disregard the exponential $e^{-\frac{|x-ry|^2}{1-r^2}}$ in \eqref{pois}  to obtain
\[ 
\int_0^{2\frac{|x|}{|y|}}\ldots\; \lesssim \, t^{2\nu}\,e^{\frac{|x|^2-|y|^2}2}\,
\int_0^{2\frac{|x|}{|y|}}
\frac{r^{\frac{m+d}2}}{\big(\ln\frac1{r}\big)^{1+\nu}}\,\frac{dr}r.
\]
This last integral
 can be estimated by $({|x|}/{|y|})^{\frac{d+m}2}/[\ln({|y|}/{|x|})]^{1+\nu}$ when\footnote{When 
$m=-d$, the integral is bounded by $1/[\ln(|y|/|x|)]^{\nu}$.} $m+d>0$.
Now, using elementary bounds on logarithms (see Lemma \ref{logs} in the Appendix) one concludes that
\Be \int_0^{2\frac{|x|}{|y|}}\ldots\,\lesssim\,
\,t^{2\nu}\,e^{\frac{|x|^2}2}\,|x|^{\frac{d+m}2}\,[\ln(|x|+e)]^{1+\nu}\,\phi(y)\,.\label{I1b}\Ee
Combining \eqref{I0}, \eqref{I1a} and \eqref{I1b} one obtains the upper bound in \eqref{pk1}.

To establish the lower bound, it suffices to integrate in the range
$0< r\leq 2\frac{|x|}{|y|}$. Note that $|y|\geq3|x|$ also implies $r\leq 2/3$, so we obtain
\[
|x-ry|\leq\,|x|+r|y|\leq 3|x|\quad\Longrightarrow\quad e^{-\frac{|x-ry|^2}{1-r^2}}\geq\,e^{-18|x|^2} \]
(using $1-r^2\geq 1/2$). The first exponential in \eqref{pois} can be handled simply by \[
\exp\Big(-\tfrac{t^2}{2\ln\frac1r}\Big)\,\geq\,\exp\big(-\tfrac{t^2}{2\ln(3/2)}\big),\quad r\in[0,\tfrac23],
\]
so all together we conclude that \Beas p_t(x,y) & \gtrsim &  t^{2\nu}\,e^{-ct^2}\,e^{-18|x|^2-\frac{|y|^2}2}\,
\int_{0}^{2\frac{|x|}{|y|}}\frac{r^{\frac{m+d}2}}{\big(\ln\frac1r\big)^{1+\nu}}\,\frac{dr}r\\
& \approx & t^{2\nu}\,e^{-ct^2}\,e^{-18|x|^2-\frac{|y|^2}2}\,\frac{(|x|/|y|)^{\frac{m+d}2}}{\Ss[\ln(|y|/|x|)]^{1+\nu}}\,\gtrsim\,
\,c_1(t,x)\,\phi(y)\,,\Eeas
for some positive function $c_1(t,x)$.
\ProofEnd

\subsection{Proof of Lemma \ref{L3.2}}\label{Lem3.2}
We split the integral defining $p_t(x,y)$, as 
\[p_t(x,y)=\int_0^{\frac12}\ldots\;+\,\int_{\frac12}^1\ldots\,\leq I_0+I_1.\]
The singularity of kernel lies in the first piece $I_0$, and in order to find 
a good estimate it will be crucial to use the formula in \eqref{poi0}\footnote{
The formula for $p_t(x,y)$ in \eqref{pois} does not make so explicit the term $|x-y|$ in the 
leading exponential.}. Suppose we are in the local region 
$|y|\leq 3\max\{|x|,1\}$.
 Then, using \eqref{expts}, we can estimate $I_0$ by
\Bea I_0 & \lesssim & t^{2\nu}\;\int_0^{\frac12}\frac{e^{-\frac{ct^2+|x-y|^2}{4s}} 
\,e^{-\frac{s|x+y|^2}4}} {s^{\frac
{d}2+1+\nu}}\,ds\nonumber\\
& \approx & \tfrac{t^{2\nu}}{(ct^2+|x-y|^2)^{\frac{d}2+\nu}}
\;\int_{\frac{ct^2+|x-y|^2}2}^\infty
e^{-u}\,e^{-\frac{(ct^2+|x-y|^2)|x+y|^2}{16u}} \; {u^{\frac
{d}2+\nu-1}}\,du,\label{I0aux}\Eea where we have changed
variables $u=(ct^2+|x-y|^2)/(4s)$. 
In the last integral we can disregard
$t$ in the exponential, and overall estimate it crudely by\[ J:=\int_0^\infty
e^{-u}\,e^{-\frac{(|x-y||x+y|)^2}{16u}} \; {u^{\frac {d}2+\nu-1}}\,du\,=\,F_{\frac d2+\nu}\big(\tfrac{|x-y||x+y|}2\big),
\]
where $F_{\sigma}(z)$ was defined in \eqref{Fnu}. As we noticed in \eqref{Knu} we can write\[
F_\sigma(z)\,=\,2^{1-\sigma} \,z^\sigma\,K_\sigma(z) \,\lesssim\,(1+z)^{\sigma-\frac12}e^{-z},\quad z>0,
\]
by the standard asymptotics of $K_\sigma$; see e.g. \cite[p. 136]{leb}.
 Thus, we obtain
\[
J\lesssim\, (1+|x-y||x+y|)^{\nu+\frac{d-1}2}\;e^{-\frac{|x-y||x+y|}2}\,.
\]
Now, $|x-y||x+y|\,\geq \,|\langle x+y,x-y\rangle|\,\geq -|x|^2+|y|^2$, so in the region $|y|\leq 3\max\{|x|,1\}$ we have
\[
J\,\lesssim \,e^{\frac{|x|^2}2}\,e^{-\frac{|y|^2}2}
\;(1+|x+y|\,|x-y|)^{\nu+\frac{d-1}2}\,\lesssim\,(1+|x|)^{2\nu+d-1}\,
\,e^{\frac{|x|^2}2}\,e^{-\frac{|y|^2}2}.
\]
Inserting this into \eqref{I0aux} we obtain the bound for the local part asserted in
the statement of the lemma. 

The estimate of $I_0$ when $|y|\geq 3\max\{|x|,1\}$ is much better, of the order $I_0\lesssim t^{2\nu}e^{-(\frac12+\ga)|y|^2}$ for some $\ga>0$; see \cite[Lemma 4.2]{GHSTV} for details. Likewise, from the arguments we already used in Lemma \ref{L3.1} one obtains a bound for $I_1\leq C_2(x)t^{2\nu}\phi(y)$ with $C_2(x)\lesssim 1/\phi(x)$.
We again refer to \cite{GHSTV} for details. This completes the proof of Lemma \ref{L3.2}. 
\ProofEnd

\section{Proof of Theorems \ref{th2} and \ref{th3}}

As we mentioned in the introduction, it suffices to give a proof of Theorem~\ref{th3}. That is,
assuming $w\in D_p(\phi)$, by which we mean
\[
\|w\|_{D_p(\phi)}:=\,\big\|w^{-\frac1p}\,\phi\,\big\|_{L^{p'}(\SR^d)}<\infty,\]
we must show that the weights $v(x)$ defined in \eqref{vx} are such that $P^*_a$ maps
$L^p(w)\to L^p(v)$ boundedly, for all $a>0$. We shall use the bound for $P^*_a$ in \eqref{P**},
namely \Bea
P^*_af(x) & \lesssim &  C_1(x)\,\Ml \big(fe^{-\frac{|y|^2}2}\big)(x) \,+\,
C_2(x)\,a^{2\nu}\,\int_{\SRd}|f(y)|\,\phi(y)\,dy\nonumber\\ & = & \,
I(x)\,+\,II(x),\label{IyII}\Eea
with $C_1(x)$ and $C_2(x)$ given explicitly in Lemma \ref{L3.2}.
We treat first the last term, which  by H\"older's inequality is bounded by
\[
II(x)\,\leq\, C_2(x)\,a^{2\nu}\,\|f\|_{L^p(w)}\,\|w\|_{D_p(\phi)}\;.
\]
So using $C_2(x)=1/\phi(x)$, we will have \Be
\|II\|_{L^p(v)} \leq\,a^{2\nu}\,\|w\|_{D_p(\phi)}\, \|f\|_{L^p(w)}\,\Big[\int_{\SRd}\frac{v(x)}{\phi(x)^p}\,dx\Big]^\frac1p
\label{IIn}\Ee with the last integral being a finite expression provided we choose
\Be
v(x)\leq\, v_1(x):=\,\frac{e^{-\frac
p2|x|^2}}{(1+|x|)^{M}}
\label{v1x}\Ee
for any $M>N_1=d+p(d+m)/2$. We remark that the weights $v(x)$ in \eqref{vx} have this property, 
at least if $N$ is sufficiently large. This is a consequence of the following elementary lemma.

\begin{lemma}
\label{lowerMloc}
If $f\in L^1_{\rm loc}(\SRd)$ satisfies $f(x)>0$, $\sae x\in\SRd$, then
\Be
\Ml f(x)\,\geq \,c_f\,(1+|x|)^{-d},\quad \forall\;x\in\SRd,\label{Mlw}\Ee
with $c_f=\frac1{|B_1|}\int_{B_1(0)}f>0$.
\end{lemma}
\Proof
Choosing $r=1+|x|$, one trivially has $B_1(0)\subset B_r(x)\cap\{|y|\leq 3\max(|x|,1)\}$, so
\eqref{Mlw} is immediate from the definition of $\Ml f(x)$ in \eqref{Ml}.

\ProofEnd
Now, using the lemma, one sees that the weights $v(x)$ defined in \eqref{vx} satisfy\[
v(x)\,\lesssim \,c'_w\,(1+|x|)^{d\al(p-1)}\,e^{-\frac
p2|x|^2}(1+|x|)^{-N},\]
hence choosing $N>N_1+d\al(p-1)$ ensures that \eqref{v1x} holds.

\

We now turn to main part $I(x)$ in \eqref{IyII}. The following 
proposition will be crucial. The result is new, and the proof is based on arguments 
due to Carleson and Jones (see \cite{CJ}).

\begin{proposition}\label{Prop_Ml1}
Let $1<p<\infty$ and $w(x)>0$ such that $w^{-\frac1{p-1}}\in L^1_{\rm loc}(\SRd)$. 
Define
\Be
V_\al(x):=\frac{\Big[\Ml(w^{-\frac{1}{p-1}})(x)\Big]^{-(p-1)\al}}{(1+|x|)^{(p-1)d\al}},\quad \mbox{ for }\al>1.\label{Val}\Ee
Then \[\Ml:L^p(w)\to L^p(V_\al)\quad \mbox{boundedly.}\] Moreover, given $\sigma<1$, if we choose $1<\al<1/\sigma$, then we also have 
$V_\al^{-\frac\sigma{p-1}}\in L^1_{\rm loc}(\SRd)$.
\end{proposition}

\Proof
Note that $\Ml(w^{-\frac1{p-1}})(x)<\infty$, $\sae x\in\SRd$, by the assumption $w^{-\frac1{p-1}}\in L^1_{\rm loc}$.
This, together with Lemma \ref{lowerMloc}, imply that
\[
0<V_\al(x)\leq c_{w}, \quad \ae x\in\SRd.\]
Now call $E_n=\{x\in\SRd\mid \Ml(w^{-\frac1{p-1}})(x)<2^n\}$, $n=0,1,2,\ldots$, and define the operators
\Be
T_ng(x):=\chi_{E_n}\,\Ml(\wp g)(x).\Ee
Note that $T_n:L^1(\wp)\to L^{1,\infty}(\SRd)$, with a uniform bound  in $n$; in fact
\Be
\Big|\Big\{T_ng(x)>\la\Big\}\Big|\leq \Big|\Big\{M(\wp g)(x)>\la\Big\}\Big|\leq \frac{c_0}\la\int_{\SRd}\wp|g|.
\label{Tn0}\Ee
Similarly, $T_n:L^\infty(\wp)\to L^\infty(\SRd)$ with $\|T_n\|\leq 2^n$, since
\Be
\big\|T_ng\big\|_\infty=\sup_{x\in E_n}\big|\Ml(\wp g)(x)\big|\,\leq 2^n\,\|g\|_\infty.\label{Tninf}\Ee
Thus, by Marcinkiewicz interpolation theorem we obtain 
\Be
\int_{E_n}|T_n(g)|^p\leq\, c_0\,2^{\frac{np}{p'}}\,\int_{\SRd}|g|^p\,\wp.\Ee
Setting $g=fw^{\frac1{p-1}}$ in the above inequality, this is the same as
\Be
\int_{E_n}|\Ml(f)|^p\leq\, c_0\,2^{n(p-1)}\,\int_{\SRd}|f|^p\,w.\label{Ml_ipol}\Ee
Now, writing $\SRd=E_0\cup[\cup_{n\geq1}E_n\setminus E_{n-1}]$, and using \[
V_\al(x)\leq c_w\mbox{ in }E_0,\mand V_\al(x)\leq 2^{-(n-1)(p-1)\al} \mbox{ if }x\notin E_{n-1},
\]
we obtain
\Beas \int_{\SRd}|\Ml f|^p\,V_\al & \leq & c_w\,\int_{E_0}|\Ml f|^p\,+\,
\sum_{n=1}^\infty2^{-(n-1)(p-1)\al}\int_{E_n}|\Ml f|^p\\  
\mbox{{\tiny (by \eqref{Ml_ipol})}}& \leq & c_wc_0\,\int_{\SRd}|f|^pw\,+\,
2^{(p-1)\al}\sum_{n=1}^\infty2^{-n(\al-1)(p-1)}\int_{\SRd}|f|^pw\,\lesssim\, \int_{\SRd}|f|^pw,\Eeas  
as we wished to show. Finally note that, when $\al\sigma<1$, the classical Kolmogorov inequality gives
\[
V_\al^{-\frac{\sigma}{p-1}}(x)\,=\,\Big[\Ml(\wp)(x)\Big]^{\al\sigma}\,(1+|x|)^{d\al\sigma}\in L^1_{\rm loc}(\SRd).\]
\ProofEnd
\BR \label{Rem_conv}
The condition $\wp\in L^1_{\rm loc}(\SRd)$ actually characterizes the property that $\Ml$ maps $L^p(w)\to L^p(V)$ for some weight $V(x)>0$. We sketch a proof of the converse in Proposition \ref{App_cov} below. \ER

Below we shall need a refinement of Proposition \ref{Prop_Ml1}, which we state now. 

\begin{proposition}\label{Prop_Ml2}
In the conditions of Proposition \ref{Prop_Ml1}, if $w(x)$ additionally satisfies
\Be
\int_{\SRd}\wp(x)\,e^{-a|x|^2}\,dx\,<\,\infty,\quad\forall\;a>0,\label{wpa}\Ee
then, for every $\sigma<1$ and $1<\al<1/\sigma$, the weight $V_\al(x)$ defined in \eqref{Val}  also satisfies
\Be
\int_{\SRd}V_\al^{-\frac{\sigma}{p-1}}(x)\,e^{-b|x|^2}\,dx\,<\,\infty,\quad\forall\;b>0.\label{Valb}\Ee\end{proposition}
\Proof
Note that \eqref{Valb} is true if we only integrate in $B_1(0)$. So, we shall consider $S_j=\{2^j\leq|x|<2^{j+1}\}$,
$j=0,1,2,\ldots$ Call $s=\al\sigma<1$ and take any $b>0$. Then
\Beas
I& = & \sum_{j=0}^\infty \int_{S_j}\big|\Ml(\wp)(x)\big|^s\,(1+|x|)^{ds}\,e^{-b|x|^2}\,dx\\
& \lesssim & \sum_{j=0}^\infty 2^{jds}\,e^{-b4^j}\,\int_{S_j}\big|\Ml\big(\wp\chi_{B_{3\cdot2^{j+1}}(0)}\big)(x)\big|^s\,dx\\
\mbox{{\tiny by Kolmogorov ineq}}& \lesssim & \sum_{j=0}^\infty 2^{jds}\,e^{-b4^j}\,|S_j|^{1-s}\,\Big\|M\big(\wp\chi_{B_{3\cdot2^{j+1}}(0)}\big)(x)\Big\|^s_{L^{1,\infty}(\SRd)}\\
& \lesssim & \sum_{j=0}^\infty 2^{jd}\,e^{-b4^j}\,\Big[\int_{B_{3\cdot2^{j+1}}(0)}\wp(x)
e^{a|x|^2}e^{-a|x|^2}\,dx\Big]^s\\
& \lesssim & \sum_{j=0}^\infty 2^{jd}\,e^{-b4^j}\,e^{6^2as4^j}\,\Big[\int_{\SRd}\wp(x)e^{-a|x|^2}\,dx\Big]^s,\Eeas
and this is finite if we choose $a<b/(36s)$.\ProofEnd

\subsection{Conclusion of the proof of Theorem \ref{th3}}

Suppose now that $w\in D_p(\phi)$, that is $\int\wp(x)\phi(x)^{p'}\,dx<\infty$, for $\phi$ as in \eqref{phi}.
This implies that $W(x)=w(x)e^{\frac p2|x|^2}$ satisfies
\Be
\int_{\SRd}W^{-\frac1{p-1}}(x)\,e^{-a|x|^2}\,dx=\,\int_{\SRd}\wp(x)\,e^{-\frac{p'}2|x|^2}\,e^{-a|x|^2}\,dx<\infty,
\label{Wa}\Ee
for all $a>0$. Now, by Proposition \ref{Prop_Ml1}, $\Ml$ maps $L^p(W)\to L^p(V_\al)$ boundedly, if we set 
\[
V_\al(x)=\frac{\Big[\Ml(W^{-\frac{1}{p-1}})(x)\Big]^{-(p-1)\al}}{(1+|x|)^{(p-1)d\al}},\quad \mbox{with}\quad \al>1.
\]
In  particular, if $f\in L^p(w)$ and we write $\tf(y)=f(y)e^{-\frac{|y|^2}2}\in L^p(W)$, we have
\[
\|\Ml\tf\|_{L^p(V_\al)}\lesssim \|\tf\|_{L^p(W)}=\|f\|_{L^p(w)}.
\]
So, recalling the value of $C_1(x)$ in Lemma \ref{L3.2}, and setting 
\Be
v(x)\leq v_0(x):=\frac{V_\al(x)\,e^{-\frac p2|x|^2}}{(1+|x|)^L}\,,
\label{v0x}\Ee
with $L\geq L_1=(2\nu+d-1)p$, we see that the term $I(x)$ in \eqref{IyII} is controlled by
\Bea
\|I(x)\|^p_{L^p(v)} &  \leq & \int_{\SRd}\frac{C_1(x)^p\,e^{-\frac p2|x|^2}}{(1+|x|)^L}\,\big|\Ml\tf(x)\big|^p\,V_\al(x)\,dx\nonumber\\
& \lesssim & \|f\|_{L^p(w)}^p.
\label{In}\Eea
So, combining \eqref{IyII}, \eqref{IIn} and \eqref{In} we have shown that $\|P^*_a f\|_{L^p(v)}\lesssim\|f\|_{L^p(w)}$,
provided \[v(x)\leq \min\{v_0(x),v_1(x)\},\]
with $v_0(x)$ and $v_1(x)$ defined in \eqref{v0x} and \eqref{v1x}.
But this inequality is clearly satisfied by the weights $v(x)$ defined in \eqref{vx}, for every $\al>1$, provided that
\[N> (p-1)d\al+\max\big\{(2\nu+d-1)p, d+p\tfrac{d+m}2\big\}\,=:\,N_0.\]
Finally notice that, since \eqref{Wa} holds, if $\sigma<1$ and $1<\al<1/\sigma$, we 
can use Proposition \ref{Prop_Ml2} to obtain\[
\int_{\SRd}v(x)^{-\frac{\sigma}{p-1}}\,\phi(x)^{p'}\,dx\,\leq\,
\int_{\SRd}V_\al(x)^{-\frac{\sigma}{p-1}}\,e^{-(1-\sigma)\frac{p'|x|^2}2}\,
(1+|x|)^\frac{\sigma N}{p-1}\,dx\,<\infty.\]
This proves that we can choose the weight $v(x)$ satisfying the inequality \eqref{Dpe}, as
asserted in the introduction. It also gives a different proof (constructive) of \cite[Theorem 1.3]{GHSTV}.
\ProofEnd

\bline {\bf Acknowledgments.} I wish to thank S. Hartzstein, T. Signes, J.L. Torrea and B. Viviani
for allowing me to use the results in our joint paper \cite{GHSTV}. This presentation is strongly influenced by many conversations with them. Special thanks also to F.J. Mart\'\i n Reyes and all the organizers of the \emph{VI CIDAMA} conference for their hospitality and nice research environment enjoyed during the workshop.

\section{Appendix}

The following elementary estimates were used in the proof of Lemma~\ref{L3.1}.

\begin{lemma}\label{logs}
Let $x,y$ be positive real numbers such that $y\geq \la\max\{x,1\}$, for some $\la>1$. Then there exist $c_\la, d_\la>0$ such that 
\Be\ln \frac yx\,\geq \, c_\la\, \frac{\ln (y+e)}{\ln(x+e)}.
\label{logineq0}\Ee
and
\Be\ln \frac yx\,\leq \, d_\la\, \ln (y+e){\ln(\frac1x+e)}.
\label{logineq1}\Ee
\end{lemma}
\Proof We first consider \eqref{logineq0}. Since $\ln y/\ln(y+e)$ is bounded above and below when $y\in[\la,\infty)$,
it suffices to prove the weaker estimate
\Be\ln \frac yx\,\geq \, c'_\la\, \frac{\ln y}{\ln(x+e)}.
\label{logineq}\Ee
Consider first the case $x\leq \sqrt\la$. Then $y\geq\la$ implies that $\sqrt y\geq x$ and hence
\[
\ln\frac yx\geq \ln\sqrt y = \tfrac12\ln y\geq \tfrac12\frac{\ln y}{\ln(x+e)},
\]
since $\ln(x+e)>1$. This proves \eqref{logineq} with $c'_\la=1/2$.
Consider now the case $x\geq \sqrt\la$, and write
\Be
\ln y = \ln\tfrac yx\,+\,\ln x.
\label{logaux}\Ee
Observe that, if $a\geq a_0>0$ and $b\geq b_0>0$, then   
\Be
\frac{a+b}{ab}=\frac1a+\frac1b\leq \frac1{a_0}+\frac1{b_0}.\label{abfact}\Ee
So, using this fact in \eqref{logaux} we see that
\[
\ln y \leq \Big(\tfrac1{\ln\la}+\tfrac1{\ln\sqrt\la}\Big)\,\ln\frac yx\,\ln x\,\leq\,\tfrac3{\ln\la}\,\ln\frac yx\,\ln (x+e),\]
which implies \eqref{logineq} with $c'_\la=(\ln\la)/3$.
The proof of \eqref{logineq1} is similar. If $x\geq1/\la$ then  
\[
\ln\tfrac yx \leq \ln (\la y)\leq 2\ln y \leq 2\ln (y+e)\,\ln(e+\tfrac1x).
\]
If $x\geq1/\la$ then \[
\ln\tfrac yx = \ln y+\ln \tfrac1x  \leq \tfrac2{\ln \la}\, \ln y\,\ln(\tfrac1x),\]
using in the last step the inequality \eqref{abfact}.
\ProofEnd

We give a proof of the converse of Proposition \ref{Prop_Ml1}, whose validity we mentioned
in Remark \ref{Rem_conv}. The arguments are similar to those in \cite{CJ}.

\begin{proposition}
\label{App_cov} Let $1<p<\infty$, and suppose that $w(x)>0$ is such that $\Ml$ maps $L^p(w)\to L^p(V)$ for some weight $V(x)>0$. Then, necessarily, $\wp\in L^1_{\rm loc}(\SRd)$.
\end{proposition}
\Proof
Given $\e>0$, we set $w_\e(x)=w(x)+\e$. Notice that every $f=\wp_\e\chi_{B_R(0)}$ belongs to $L^p(w)$ since
\[
\int_{\SRd}|f|^pw\leq \int_{B_R(0)}w^{-p'}_\e\,w_\e\,\leq\,|B_R(0)|\,\e^{-(p'-1)}\,<\,\infty.\]
Call $S_0=B_1(0)$ and $S_j=\{2^{j-1}\leq|x|<2^j\}$, $j=1,2,\ldots$, and consider the $L^p(w)$-functions 
$f_j=\wp_\e\chi_{S_j}$. Then, it is easy to verify from the definition of $\Ml$ that
\[
\Ml f_j(x) \gtrsim 2^{-jd}\int_{S_j} \wp_\e,\quad \mbox{if }x\in S_j,\]
for each $j=0,1,2,\ldots$ Indeed, it suffices to average over a ball $B_r(x)$ of radius $r=2^j+|x|$,
and observe that $S_j\subset B_r(x)\cap\{|y|\leq3\max(|x|,1)\}$ when $x\in S_j$.
Thus, the assumed boundedness of $\Ml$ gives\[
V(S_j)^\frac1p\,2^{-jd}\int_{S_j} \wp_\e\lesssim \big\|\Ml f_j\big\|_{L^p(V)}\,\leq\, C\, \|f_j\|_{L^p(w)}\,\leq\,
C\Big[\int_{S_j} w_\e^{-(p'-1)}\Big]^{\frac1p},\]
with a constant $C$ independent of $\e$. Since both sides are finite and $p'-1=1/(p-1)$, we conclude that
\[
\Big[\int_{S_j} \wp_\e\Big]^{1/{p'}}\,\lesssim \, C\,2^{jd}/V(S_j)^\frac1p,
\]
so letting $\e\to 0$ we obtain that $\wp\in L^1(S_j)$ for each $j=0,1,2,\ldots$\ProofEnd

\bibliographystyle{plain}

\end{document}